\documentclass{amsart}

\def\kinf{\mbox{$K_{\infty}$}}
\def\univk{\mbox{$\mathcal{U}^k$}}
\def\gnalpha{\mbox{$G(n,n^{-\alpha})$}}
\def\foxk{\mbox{$\preceq^{k}$}}

\def\exk{\mbox{$\preceq_{\infty\omega}^k$}}

\def\nxk{\mbox{$\preceq^{k,n}$}}
\def\n1xk{\mbox{$\preceq^{k,n+1}$}}
\def\dn1xk{\mbox{$\not\preceq^{k,n+1}$}}

\def\*{\makebox[10mm]{}}

\def\Co-NP{\mbox{\rm Co-NP}}

\def\as{\mbox{\rm a.s.}}

\newcommand{\cn}[1]{\mbox{$\chi^{\ast}(#1)$}}
\newcommand{\ffoe}[1]{\mbox{$L^{#1}(\exists)$}}
\newcommand{\ki}[1]{\mbox{$K^{\infty}_{#1}$}}
\newcommand{\qr}[1]{\mbox{${\rm qr}(#1)$}}
\newcommand{\edges}[1]{\mbox{${\rm Edges}(#1)$}}

\newcommand{\safesub}[2]{\mbox{$#1 \leq_s #2$}}
\newcommand{\rigidsub}[2]{\mbox{$#1 \leq_i #2$}}
\newcommand{\speciall}[3]{\mbox{$#2 \leq_{#1}^{\otimes} #3$}}
\newcommand{\freejoin}[3]{\mbox{$#1 \otimes_{#3} #2$}}
\newcommand{\speciallns}[2]{\mbox{$#1 \leq^{\otimes} #2$}}
\newcommand{\cl}[4]{\mbox{${\rm cl}^{#1,#2}(#3,#4)$}}

\def\leqi{\mbox{$\leq_{i}$}}

\def\ekgame{\mbox{$\exists^k$-game}}
\def\ekpgame{\mbox{$\exists^{k'}$-game}}

\def\finmodalpha{\mbox{$K^{\alpha}_{\infty}$}}

\def\fok{\mbox{$L^k$}}
\def\foke{\mbox{$L^k(\exists)$}}
\def\lomega{\mbox{$L^{\omega}_{\infty\omega}$}}
\def\lk{\mbox{$L^k_{\infty\omega}$}}
\def\lke{\mbox{$L^k_{\infty\omega}(\exists)$}}

\def\linf{\mbox{$L_{\infty\omega}$}}

\def\orr{\vee}

\def\cnk{\mbox{$\bigwedge L^{k}(\exists)$}}
\def\dsk{\mbox{$\bigvee L^{k}(\exists)$}}

\newcommand{\dke}[1]{\mbox{$\bigvee L^{#1}(\exists)$}}

\newtheorem{theorem}{Theorem}[section]
\newtheorem{lemma}[theorem]{Lemma}
\newtheorem{sublemma}[theorem]{Sublemma}
\newtheorem{corollary}[theorem]{Corollary}
\newtheorem{proposition}[theorem]{Proposition}

\theoremstyle{definition}
\newtheorem{definition}[theorem]{Definition}

\theoremstyle{remark}

\numberwithin{equation}{section}

\begin{document}

\title{$k$-Universal Finite Graphs}

\author{Eric Rosen}
\address{Department of Computer Science, Technion--Israel Institute of
Technology, Technion City, Haifa 32000, Israel} 
\email{erosen@csa.cs.technion.ac.il}
\thanks{The first author was supported in part by a postdoctoral fellowship
from 
the Israel Council for Higher Education and by a research grant of the German
Israeli Foundation (GIF)}

\author{Saharon Shelah}
\address{Institute of Mathematics, The Hebrew University of Jerusalem, 91904
Jerusalem, Israel {\rm and} Department of Mathematics, Rutgers University, New
Brunswick, New Jersey 08854}
\email{shlhetal@sunset.huji.ac.il}
\thanks{The second author's research was supported by the United States-Israel
Binational Science Foundation and by DIMACS. Publication No.~611.}

\author{Scott Weinstein}
\address{Department of Philosophy, University of Pennsylvania, Philadelphia,
Pennsylvania 19104}
\email{weinstein@cis.upenn.edu}
\thanks{The third author was supported in part by NSF CCR-9403447.
We would like to thank John Baldwin and an anonymous referee for useful
comments on an earlier draft of this paper.}

\date{April 9, 1996}

\begin{abstract}
This paper investigates 
the class of $k$-universal finite graphs,
a local analog of the class of universal graphs,
which arises naturally in the study of finite variable logics. 
The main results of the paper, which are due to Shelah, establish that the
class of $k$-universal graphs is not definable by an infinite disjunction of
first-order existential
sentences with a finite number of variables and that there exist
$k$-universal graphs with no $k$-extendible induced subgraphs.
\end{abstract}

\maketitle
\section{Introduction}
This paper continues the investigation of the existential fragment of \lomega\
from the point of view of finite model theory initiated in \cite{rnw} and
\cite{Rosen-diss}. 
In particular, we further study an analog of universal structures, namely,
$k$-universal structures, which arise
naturally in the context of finite variable logics.
The main results of this paper, Theorems \ref{main-thm} and
\ref{univnotextend-thm},
which are due to Shelah, 
apply techniques from the theory of sparse
random graphs as developed in \cite{shelah-spencer} and \cite{baldwin-shelah}
to 
answer some
questions about $k$-universal structures left open in these earlier works. 
In order to make the current paper more or less self-contained, we recall some
notions and notations
from the papers cited above, which may be consulted for further
background and references.

We restrict our attention to languages which contain only relation symbols. 
We let \fok\ denote the fragment of first-order logic consisting of those
formulas 
all of whose variables both free and bound are among $x_1, \ldots, x_k,$
and similarly, \lk\ is the $k$-variable fragment of the infinitary language
\linf. 
We let \foke\ denote the collection of existential formulas of \fok,
that is, those
formulas 
obtained by closing the set of atomic 
formulas and negated atomic formulas of \fok\ under the operations of 
conjunction, disjunction, and existential quantification, and we let \lke\ be
the existential fragment of \lk.
The fragments \cnk\ and \dsk\ of \lke\ consist of the countable
conjunctions and the countable disjunctions of formulas of \foke\
respectively. We write \qr{\theta}\ for the quantifier rank of the formula
$\theta,$ which is defined as usual. 
\begin{definition}\label{compatible-def}
Let $A$ and $B$ be structures of the same relational signature. $A \foxk B$ ($A
\nxk B$) ($A \exk B$), if and only if, 
for all $\theta\in \foke$ (with $\qr{\theta} \leq n$) (for all $\theta \in
\lke$), if $A \models \theta,$
then $B \models \theta.$ 
\end{definition}
These relations may be usefully characterized in terms
of the following non-alternating, local variants of the 
Ehrenfeucht-Fraisse 
game. 
The {\em $n$-round}, \ekgame\ from $A$ to 
$B$ is played between two players, Spoiler and Duplicator, with $k$
pairs of pebbles, $(\alpha_{1}, \beta_{1}), \ldots ,$ 
$(\alpha_{k}, \beta_{k})$.
The Spoiler begins each round by choosing a pebble $\alpha_{i}$
 that may or may not be in play and placing it on an element of $A.$
The Duplicator then plays $\beta_i$ onto an element of $B.$
The Spoiler wins the game if after any round $m \leq n$ the function $f$ from $A$ to
$B,$ which sends the element pebbled by $\alpha_{i}$ to the 
element pebbled by
$\beta_i$ is not a partial
isomorphism; otherwise, the Duplicator wins the game.
The {\em eternal} \ekgame\ is an infinite version of the $n$-round game in
which 
the 
play continues through a sequence of rounds of order type $\omega.$ The Spoiler
wins the game, if and only if, he wins at the $n^{\rm th}$-round for some $n
\in \omega$ as above; otherwise, the Duplicator wins. 
The following proposition provides the link between the \ekgame\ and logical
definability.
\begin{proposition}[\cite{KV-datalog}] \label{ekgame-prop}
\begin{enumerate}
\item For all structures $A$ and $B,$ the following conditions are equivalent.
\begin{enumerate}
\item $A \nxk B.$ 
\item The Duplicator has a winning strategy for the
$n$-round \ekgame\ from $A$ to $B.$ 
\end{enumerate}
\item For all structures $A$ and $B,$ the following conditions are equivalent.
\begin{enumerate}
\item $A \exk B.$
\item The Duplicator has a winning strategy for the 
eternal \ekgame\ from $A$ to $B.$
\end{enumerate}
\item For all structures $A$ and finite structures $B,$ the following
conditions are equivalent. 
\begin{enumerate}
\item $A \exk B.$
\item $A \foxk B.$
\end{enumerate}
\end{enumerate}
\end{proposition}

In this paper, we will focus our attention on the class of finite simple
graphs, that 
is, finite structures with one binary relation which is irreflexive and
symmetric. We will use the term graph to refer to such structures. 
In general, we let $A, B, \ldots$ refer both to graphs and to their underlying
vertex sets and we let $|A|$ denote the cardinality of $A.$ We use
$E$ for the edge relation of a graph. \edges{A} is the edge set of the graph
$A,$ that is,  $\edges{A} = \{ \{a,b\}\subseteq A : E(a,b) \}.$ 

\section{$k$-Universal Graphs: Definability and Structure}

We say that a graph $G$ is {\em k-universal}, if and only if, for all graphs
$H, H \foxk G.$  By Proposition \ref{ekgame-prop}, this is equivalent to $G$
satisfying every sentence of \lke\ which is satisfied by some (possibly
infinite) 
graph. 
We say that a graph $G$ is {\em k-extendible}, if
and only if, $k \leq |G|$ and 
for each $1 \leq l \leq k$
\begin{equation*}
\begin{split}
G \models \forall x_1 \ldots \forall x_{k-1}\exists x_k 
&( \bigwedge_{1 \leq i
< j \leq 
k-1} x_i \neq x_j \rightarrow \\
&(\bigwedge_{1 \leq i \leq 
k-1} x_i \neq x_k \wedge
\bigwedge_{1 \leq i < l} E(x_i,x_k) \wedge \bigwedge_{l \leq i < k}
\neg E(x_i,x_k))). 
\end{split}
\end{equation*}
It is easy to verify, by applying Proposition \ref{ekgame-prop}, that
every $k$-extendible graph is $k$-universal. 
The class of $k$-extendible graphs plays an important role in
the study of $0-1$ laws for certain infinitary logics and logics with
fixed point operators (see \cite{KV-IC}). 
Indeed, the existence of $k$-universal
finite graphs follows immediately from the fact that for every $k,$
the random graph $G = G(n,p)$ with constant edge probability $0 < p < 1$ is
almost surely $k$-extendible (see, for example, \cite{Bollobas-text}).

Let \univk\ be the class of $k$-universal graphs and let 
$$\Xi^k = \{ \theta \in
\foke : \exists G(G\ {\rm is\ a\ graph\ and\ } G \models \theta) \}.$$
Note that
for all graphs $G, G \in \univk,$ if and only if, $G \models \bigwedge \Xi^k.$
Thus, \univk\ is definable in \cnk\ over the class of graphs.   
In \cite{rnw}, we established via an explicit construction
that for all $2 \leq k,$ \univk\ is not definable
in \dsk. The following theorem 
significantly strengthens this result for large
enough $k$; its proof involves a probabilistic construction employing
techniques 
from the theory of sparse random graphs.  
\begin{theorem}\label{main-thm}
For all $k \geq 7$ and $k'\in \omega,$
\univk\ is not definable in \dke{k'} over
the class of graphs.
\end{theorem}

We call a class of structures $\mathcal{C}$ finitely based, if and only if, there
is a finite set of structures $\{A_1, \ldots , A_n \} \subseteq \mathcal{C}$ such
that for every structure $B \in \mathcal{C}, A_i \subseteq B$ for some $ 1 \leq i
\leq n.$ We obtain the following result as a corollary to the proof of Theorem
\ref{main-thm}. 
\begin{corollary}\label{finbase-cor}
For all $k \geq 7,$ 
\begin{enumerate}
\item \univk\ is not finitely based, and
\item the class of $k$-extendible graphs is not finitely based.
\end{enumerate} 
\end{corollary} 

In \cite{rnw}, we observed that for all $k,$ \univk\ is decidable in
deterministic polynomial time. The following theorem gives a stronger
``descriptive complexity'' result. 
\begin{theorem}\label{univinlfp-thm}
For all $k,$ \univk\ is definable in least fixed point logic.
\end{theorem}

It is clear that if $G$ is $k$-extendible and $G \subseteq H,$ then $H$ is
$k$-universal. The question naturally arises whether 
there are $k$-universal graphs which contain no $k$-extendible subgraph.
The following theorem answers this question
affirmatively.
\begin{theorem}\label{univnotextend-thm}
For each $k \geq 4,$ there is a graph $G$ such that 
\begin{enumerate}
\item $G$ is $k$-universal, and
\item $\forall H \subseteq G, H$ is not $k$-extendible.
\end{enumerate}
\end{theorem}

The next theorem is a strengthening of the first part of Corollary
\ref{finbase-cor}. The proof of this theorem expands on the construction
developed to prove Theorem \ref{univnotextend-thm}. We say a graph $G$ is a
{\em minimal} $k$-universal graph just in case $G$ is $k$-universal and
contains no proper induced subgraph which is $k$-universal.

\begin{theorem} \label{kineqminkuniv-thm}
For all $k \geq 6$, there is an infinite set of pairwise
$L^k$-inequivalent minimal $k$-universal graphs.
\end{theorem}
 
We proceed to prove the above results.
Theorem \ref{main-thm} is an immediate corollary of the following lemma
which
is due to Shelah.

\begin{lemma}\label{central-lemma}
For all $k \geq 7$ and $k' \in \omega,$ 
there is a 
graph
$N$ such that
\begin{enumerate}
\item $N$ is $k$-extendible and
\item for every $\theta \in \ffoe{k'},$ if $N \models \theta,$ then there is a
structure $M$ such that $M \models \theta$ and $M$ is not $k$-universal.
\end{enumerate}
\end{lemma}

We approach the proof of Lemma \ref{central-lemma} through a sequence of
sublemmas. We first introduce some graph-theoretic concepts 
which play a central role in the argument. 
\begin{definition}\label{coloring-number-def}
Let $A$ be a finite graph.
\begin{enumerate}
\item We say
$\overline{a}=\langle a_1, \ldots , a_n\rangle$ is a {\em t-witness} for $A$,
if and only 
if, $\overline{a}$ is an injective enumeration of $A$ and
for each $i 
\leq n, |\{ j < i : E(a_j,a_i) \}| \leq t.$
\item $\cn{A} = $\ the least $t$ such that there is a $t$-witness for $A.$
(\cn{A} is the {\em coloring number} of $A.$)
\item $\ki{t} = \{ A : \cn{A} \leq t \}.$
\item \speciall{t}{A}{B}, if and only if, $A \subseteq B, 
B \in \ki{t}$ and every $t$-witness for
$A$ can be extended to a $t$-witness for $B,$ that is, if $\overline{a}$ is a
$t$-witness for $A,$ then there is a $\overline{b}$ such that
$\overline{a}\overline{b}$ is a $t$-witness for $B.$
\end{enumerate}
\end{definition}
The coloring number was
introduced and extensively studied in \cite{erdos-hajnal-66}.
The following sublemma states a free amalgamation property of
$\speciall{t}{}{}.$
\begin{definition}
Let $A$ and $B$ be finite graphs. 
\begin{enumerate}
\item $A$ is {\em compatible} with $B$, if and only if, the subgraph of $A$
induced by $A \cap B$ is identical to the subgraph of $B$ induced by $A \cap
B.$ 
\item Suppose $A$ is compatible with $B$ and let $C$ be the subgraph of $A$
induced by $A \cap B.$
The {\em free join} of $A$ and $B$ over $C,$
denoted by $\freejoin{A}{B}{C},$ is the graph whose vertex set is $A \cup B$
and whose edge set is $\edges{A} \cup \edges{B}.$
\end{enumerate}
\end{definition}
\begin{sublemma}\label{freeamalgam-sublemma}
Suppose $A,B \in \ki{t},$
$A$ is compatible with $B,$ $C$ is the subgraph of $A$
induced by $A \cap B,$ $\speciall{t}{C}{A},$ and
$\speciall{t}{C}{B}.$ Then, $\freejoin{A}{B}{C} \in \ki{t},$
$\speciall{t}{A}{\freejoin{A}{B}{C}},$ and
$\speciall{t}{B}{\freejoin{A}{B}{C}}.$ 
\end{sublemma}
\begin{proof}
The sublemma follows immediately from the definitions.
\end{proof}

The next sublemma establishes a lower bound on \cn{G}\ when $G$ is
$k$-universal. For the proof of the sublemma we extend the definition of 
$k$-universality to apply also to tuples.
We also introduce a 
refinement of the concept that will be used in the 
proof of Theorem 2.  An $m$-tuple $\overline{a}
= (a_1, \ldots , a_m )$ is {\em proper} iff
for all $i < j \leq m, a_i \neq a_j$.  For all models $A$ and 
$B$, and $j$-tuples $\overline{a} \subseteq A, \overline{b} 
\subseteq B$, we write $(A, \overline{a})\foxk (B, \overline{b}) 
((A, \overline{a})\nxk (B, \overline{b}))$ iff for all formulas 
$\theta(\overline{x}) \in \foke$ (with $qr(\theta) \leq n$),
with $j$ free variables,
if $A \models \theta[\overline{a}]$, then $B \models \theta[\overline{b}]$.

\begin{definition}
For $j \leq k$, a proper $j$-tuple $\overline{a} \subseteq A$ is
{\em $k$-universal in $A$} ({\em $k,n$-universal in $A$}) iff for all $B$, and
proper $j$-tuples  
$\overline{b} \subseteq B$ such that the partial function 
$f(x)$ from $A$ to $B$ that maps $a_i$ to $b_i$ is a partial 
isomorphism,
$(B,\overline{b}) \foxk (A, \overline{a})$
($(B,\overline{b}) \nxk (A, \overline{a})$).  The {\em rank} of
$\overline{a} \subseteq A$ is $\omega$ if it is $k$-universal, and
the greatest $n$ such that it is $k,n$-universal, otherwise.
\end{definition}
\begin{sublemma}\label{cnuniv-sublemma}
If $\cn{G} < 2^{k-2},$
then $G$ is not $k$-universal.
\end{sublemma}
\begin{proof} 
Suppose $\cn{G} < 2^{k-2},$ and, for
{\em reductio}, that $G$ 
is $k$-universal. 
Suppose $G = \{ a_i : i < n \},$ and
let $$I = \{ \langle i_1, \ldots , i_k \rangle : i_1 < i_2
< \ldots 
< i_k < n\ {\rm and}\ \langle a_{i_1}, \ldots, a_{i_k}\rangle\ {\rm is}\ k{\rm-universal\ in}\ G \}.$$ Since $G$ is $k$-universal, it follows that $I \neq
\emptyset.$ Let $\langle i_1, \ldots , i_k \rangle \in I$ with $i_k$
maximal. Let $w = \{j < 
i_k : E(a_j,a_{i_k}) \},$ and for each $j \in w,$ let $u_j = \{ l : l \in
\{1, \ldots , k-1 \}\ {\rm and}\ E(a_j,a_{i_l}) \}.$ Choose $l^{\ast} \in \{1, \ldots ,
k-1 \}.$ As $ |w| < 2^{k-2},$ there is $u \subseteq \{1, \ldots, k-1 \}
- \{ l^{\ast} \}$ such that for every $j \in w, u \neq u_j - \{ l^{\ast} \}.$ 

Now, let $H$ be a $k$-extendible graph with edge relation $E'.$ Since $\langle
a_{i_1}, 
\ldots, 
a_{i_k}\rangle$ is $k$-universal 
in $G,$ we may choose $b_1, \ldots , b_k \in H$ such that the
Duplicator 
has a winning strategy for the \ekgame\ played from $H$ to $G$ with the $j^{\rm
th}$ pair of pebbles placed on $b_j$ and $a_{i_j}.$ We show that, in fact, the
Spoiler can force a win from this position, which yields the desired
contradiction. The Spoiler picks up the pebble resting on $b_{l^{\ast}}$ and places
it on a point $b \in H - \{ b_1 , \ldots , b_k \}$ such that $E'(b,b_k)$ and
$E'(b,b_l)$ for each $l \in u$ 
while $\neg E'(b,b_l)$ for each $l \in \{1, \ldots , k-1\} - (u \cup \{ l^{\ast} \}).$
In order to successfully answer the Spoiler's move, the Duplicator must move the
pebble now resting on $a_{i_{l^{\ast}}}$ and place it on a point $a_m \in G$ such
that $E(a_m,a_{i_k})$ and $a_m \neq a_{i_k}.$ In order to achieve this, she must
choose $a_m$ so that either $i_k < m$ or $m \in w.$ But in the first case we
would have that the position $\langle \ldots, \langle b_j,a_{i_j}\rangle , \ldots ,
\langle b_k,a_{i_k}\rangle,\langle b,a_m\rangle : j \neq l^{\ast}\rangle$ is a winning position for the
Duplicator in the \ekgame\ from $H$ to $G.$ This implies that $\langle \ldots ,
a_{i_j} , \ldots , a_{i_k} , a_m : j \neq l^{\ast}\rangle$ is $k$-universal
in $A.$ But then, since $i_k < m,$ we have 
$\langle \ldots, i_j , \ldots, i_k , m : j
\neq l^{\ast} \rangle \in I.$ But, this contradicts the choice of $i_k$ to be maximal with
this property. Therefore, it suffices to show that $m \not\in w.$ But this
follows immediately from the fact that $m < i_k$ and the
construction of $u.$ 
\end{proof}

The next sublemmas deal with the theory of the random graph $G=\gnalpha$,
$\alpha$ an irrational between $0$ and $1,$ as 
developed in \cite{shelah-spencer} (see also 
\cite{baldwin-shelah} for connections with model theory).
We say a property holds almost surely (abbreviated \as) in \gnalpha, if and
only if, its probability approaches $1$ as $n$ increases. Shelah and Spencer
showed (see \cite{shelah-spencer}) that for any first-order property $\theta$
and any irrational $\alpha$ between $0$ and $1,$ either $\theta$ holds \as\ in
\gnalpha\ or $\neg\theta$ holds \as\ in \gnalpha. For each such $\alpha,$ we
let $T^{\alpha} = \{ \theta : \theta\ {\rm holds\ \as\ in\ } \gnalpha
\}$ and we let \finmodalpha\ be the set of finite graphs each of which is
embeddable in every model of $T^{\alpha}.$ We will suppress the superscripts on
these notations, when no confusion is likely to result; in general, we will use
notations which leave reference to a particular $\alpha$ implicit, as in the
following definition. 
\begin{definition}[\cite{shelah-spencer}]\label{sparse-def}
Let $G$ and $H$ be graphs with $G \subseteq H,$ and let $\alpha$ be a fixed
irrational between $0$ and $1.$
\begin{enumerate}
\item $(G,H)$ is {\em sparse}, if and only if, $|\edges{H} - \edges{G}|/|H - G|
< 1/\alpha.$
\item $(G,H)$ is {\em dense}, if and only if, $|\edges{H} - \edges{G}|/|H - G|
> 1/\alpha.$
\item \safesub{G}{H}, if and only if, for every $I,$ if $G \subset I \subseteq
H,$ then $(G,I)$ is sparse.
\item \rigidsub{G}{H}, if and only if, for every $I,$ if $G \subseteq I \subset
H,$ then $(I,H)$ is dense.
\end{enumerate}
We say $G$ is {\em sparse} ({\em dense}), if and only if, $(\emptyset,G)$ is
sparse (dense).
\end{definition}

Note that since $\alpha$ is irrational every $(G,H)$ as above is either sparse
or dense.

\begin{sublemma}\label{finmod-sparse-sublemma}
If $G \in \kinf,$ then $\safesub{\emptyset}{G}.$
\end{sublemma}
\begin{proof} 
The reader may find a proof of this sublemma in
\cite{spencer-90}.
\end{proof} 

\begin{sublemma}\label{finmod-cn-sublemma}
If $\alpha$ is irrational and $1/(k+1) < \alpha < 1,$ then
\begin{enumerate}
\item  $\kinf \subseteq \ki{(2k+1)}$ and
\item if \safesub{A}{B}, then \speciall{2k+1}{A}{B}.
\end{enumerate}
\end{sublemma}
\begin{proof} 
1. By Sublemma \ref{finmod-sparse-sublemma}, it suffices to show
that if $\safesub{\emptyset}{G},$
then $G \in \ki{2k+1}.$ 
So suppose $\safesub{\emptyset}{G}.$ We inductively define a $2k+1$-witness for
$G$ proceeding from the top down. Since $G$ is sparse,
$|\edges{G}|/|G| < k+1,$ from which it follows immediately that there is a
point $a \in G$ whose degree is $< 2k+2.$ We let $a = a_{|G|}$ be the last
element of our $2k+1$-witness for $G.$ Now, since
$\safesub{\emptyset}{G},$ $G' = G - \{a\}$ is sparse, so we may find an $a' \in
G'$ whose degree (in $G'$) is $< 2k+2$ as before. We let $a' = a_{|G|-1}$ be
the next to last element of our $2k+1$-witness for $G.$ Proceeding in this way,
we may complete the construction of a $2k+1$-witness for $G.$ 

2. Suppose \safesub{A}{B} and suppose $\overline{a}$ is a $2k+1$-witness for
$A.$ Just as above we may inductively construct an enumeration $\overline{b}$
of $B-A$ so that $\overline{a}\overline{b}$ is a $2k+1$-witness for $B.$
\end{proof}

The following closure operator plays an important role in the proof of Lemma
\ref{central-lemma}.

\begin{definition}\label{closure-def}
We define for graphs $G, H$ with $G \subseteq H$ and  natural numbers $l,$ a
closure operator $\cl{l}{m}{G}{H}$\ by recursion on $m.$
\begin{enumerate}
\item $\cl{l}{0}{G}{H} = G;$
\item $\cl{l}{m+1}{G}{H} = \bigcup \{ B : B \subseteq H\ {\rm and}\ |B| \leq l\
{\rm and}\ B \cap 
\cl{l}{m}{G}{H} \leqi B \}.$ 
\end{enumerate}
We let $\cl{l}{\infty}{G}{H} = \bigcup_{m \in \omega}\cl{l}{m}{G}{H}.$
We say that $H$ is $l$-{\em small}, if and
only if, there is a $G \subseteq H$ such that $|G| \leq l$ and 
$\cl{l}{\infty}{G}{H} = H.$
\end{definition}

The following lemma gives the crucial property of closures we will
exploit -- for a fixed $l$ there is almost surely in \gnalpha\ a uniform
bound on the cardinality of the closure of a set of size at most $l.$

\begin{sublemma}\label{closure-sublemma}
For every $l$ there
is an 
$l^{\ast}$ such that \as\
for every $A \subseteq G (= \gnalpha)$, if $|A| \leq l,$ then 
$|\cl{l}{\infty}{A}{G}| \leq l^{\ast}.$ 
\end{sublemma}
\begin{proof} 
Note that if \rigidsub{B}{B'} and $B \subseteq C \subseteq B',$
then \rigidsub{C}{B'}. It follows that we may represent \cl{l}{\infty}{A}{G} as
$A \cup \bigcup_{i < i^{\ast}} B_i$ where $|B_i| \leq l$ and  
$(A \cup \bigcup_{j < i} B_j) \cap B_i \leq_i B_i.$ 
Moreover, we may suppose,
without loss of generality, that this last extension is strict, for otherwise
$B_i$ could be omitted from the representation. Next we argue that there is an
$m$ (depending on $l$) which \as\ uniformly bounds $i^{\ast},$ that is, there
is 
an $m$ such that 
\begin{quote}
$(\dagger)$ \as\ in $G = \gnalpha$ for all $A \subseteq G, |A| \leq l,$
there is an $i^{\ast} \leq m$ such that
 \cl{l}{\infty}{A}{G} may be represented as
$A \cup \bigcup_{i < i^{\ast}} B_i$ where $|B_i| \leq l$ and  
$\rigidsub{(A \cup \bigcup_{j < i} B_j) \cap B_i}{B_i}.$ 
\end{quote}
The sublemma follows 
immediately from this, for then $l^{\ast} = m\cdot l$ is an \as\ uniform bound
on $|\cl{l}{\infty}{A}{G}|.$ 

Let
\begin{equation*}
\begin{split}
\varepsilon = {\rm Min}(\{& (\alpha \cdot |\edges{B}-\edges{C}|)
-(|B-C|): \\
&B \subseteq G, |B| \leq l, \rigidsub{A \cap B}{B}, A \cap B \subseteq C \subset
B  \}).
\end{split}
\end{equation*}
It follows from the definition of $\leq_i$ that
$\varepsilon > 
0.$
Let $m = 1 + l/\varepsilon.$ We claim that $m$ satisfies condition $(\dagger).$
Let $$w_i = |A \cup \bigcup_{j < i} B_j| - \alpha \cdot 
|\edges{A \cup \bigcup_{j < i} B_j}|.$$
Then, by hypothesis, $w_0 \leq |A| \leq l.$ Moreover, $w_{i+1} \leq (w_i -
\varepsilon).$ To see this, let $C = B_i \cap (A \cup \bigcup_{j < i} B_j)).$
Then, $A \cap B_i \subseteq C \subset B_i.$ Hence,
$w_{i+1} = |(A \cup \bigcup_{j < i} B_j) \cup B_i| - \alpha \cdot |\edges{(A
\cup \bigcup_{j < i} B_j) \cup B_i}| \leq (|A \cup \bigcup_{j < i} B_j| + |B_i
- C|) - \alpha \cdot (|\edges{A
\cup \bigcup_{j < i} B_j}| +(|\edges{B_i}| - |\edges{C}|)) \leq (w_i -
\varepsilon).$ 
It follows, by induction, that $w_i \leq l - i\cdot\varepsilon.$ Therefore,
if $i > l/\varepsilon,$ then $w_i < 0.$ So, by Sublemma
\ref{finmod-sparse-sublemma},
if $i^{\ast} \geq m,$ 
then $\cl{l}{\infty}{A}{G} = 
A \cup \bigcup_{i < i^{\ast}}B_i \not\in \kinf.$ 
Therefore, \as\ 
$i^{\ast} < m.$
\end{proof}

For the purposes of the next sublemma and beyond, we introduce the
following notational convention: 
we write \speciallns{A}{B}\ for \speciall{t}{A}{B}, when $t = 2^{k-2}-1.$
 
\begin{sublemma}\label{cn-closure-sublemma}
If $\alpha$ is irrational, $1/(k+1) < \alpha < 1, k \geq 7$ and $k+1 < k'$ then
the following condition holds \as\ in $G = \gnalpha.$ 
For all $a_1,
\ldots, a_{k'} \in G,$ if $A = \cl{k'}{\infty}{\{a_1,\ldots,a_{k'-1}\}}{G}$
and $B  = \cl{k'}{\infty}{\{a_1,\ldots,a_{k'}\}}{G},$ then 
\begin{enumerate}
\item $B \in \kinf$ and 
\item \speciallns{A}{B}.
\end{enumerate}
\end{sublemma}
\begin{proof} 
1. This is an immediate consequence of the preceding Sublemma. By
the first-order 0-1 law for \gnalpha, given any fixed bound $l^{\ast},$ \as\ 
for all $A \subseteq G,$ if $|A| \leq l^{\ast},$ then $A \in \kinf.$ \\
2. First observe that our closure operator is monotone in $\subseteq,$ hence $A
\subseteq B$ and also, by the definition of the closure operator, that for no
$C \subseteq B, C \not\subseteq A,
|C| \leq k'$ do we have $\rigidsub{A \cap C}{C}.$ We argue that
\speciallns{A}{B} as follows. Suppose $\overline{a} = \langle a_1, \ldots, a_{|A|}\rangle$
is a 
$2^{k-2}-1$-witness for $A,$ and let $\overline{b} = \langle b_1, \ldots, b_{|B|}\rangle$ be
a 
$2k+1$-witness for $B.$ The latter exists by Sublemma 
\ref{finmod-cn-sublemma} since $B \in
\kinf.$ Now, for every $b \in B-A, |\{a\in A : E(a,b)\}| \leq k,$ for otherwise
we could find a set $C \subseteq B, C \not\subseteq A, |C| = k+2,$ such that
\rigidsub{A \cap 
C}{C}. Let $w =\{ i : 1 \leq i \leq |B|\ {\rm and}\ b_i \not\in A \},$ and let
$\overline{b'} = \langle b_i : i \in w\rangle$ be the restriction of $\overline{b}$ to an
enumeration of $B-A.$ By hypothesis, $k \geq 7,$ so $(2k+1)+k \leq
2^{k-2}-1;$ hence, we may 
conclude that $\overline{a}\overline{b'}$ is a $2^{k-2}-1$-witness for $B.$
\end{proof}

\begin{sublemma}\label{alpha-kextendible-sublemma}
If  $0 < \alpha < 1/k,$ then \gnalpha is \as\ $k$-extendible.
\end{sublemma}
\begin{proof} 
The reader may find a proof of this sublemma in 
\cite{mcarthur-diss}. 
\end{proof}

We are now in a position to proceed to the proof of Lemma \ref{central-lemma}.

\begin{proof}[Proof of Lemma \ref{central-lemma}] 
Let $k \geq 7$ and, without
loss 
of 
generality, let $k' > k+1.$ 
Fix $\alpha$ to be an irrational
number between $1/(k+1)$ and $1/k.$ It then follows from Sublemmas
\ref{cn-closure-sublemma} and
\ref{alpha-kextendible-sublemma}  that there is a finite graph $N$ such that
\begin{enumerate}
\item[(N1)] $N$ is $k$-extendible;
\item[(N2)] for all $a_1,
\ldots, a_{k'} \in N,$ if $A = \cl{k'}{\infty}{\{a_1,\ldots,a_{k'-1}\}}{N}$
and \\
$B  = \cl{k'}{\infty}{\{a_1,\ldots,a_{k'}\}}{N},$ then $B \in \kinf$ and 
\speciallns{A}{B}.
\end{enumerate}

To complete the proof we must construct for each $\theta \in \ffoe{k'},$ a
graph $M$ such that $M$ is not $k$-universal and if $N \models \theta,$ then $M
\models \theta.$ By Sublemma \ref{cnuniv-sublemma} and Proposition
\ref{ekgame-prop}, it suffices to construct for each $d \in \omega$ a graph $M$
such that
\begin{enumerate}
\item[(M1)] $\cn{M} < 2^{k-2},$ and 
\item[(M2)] the Duplicator has a winning strategy for the $d$-move \ekpgame\
from $N$ to 
$M.$ 
\end{enumerate}

We proceed to construct a structure $M$ that satisfies conditions (M1) and
(M2). 
We first define chains of structures 
 $\langle M_i : i \leq d+1\rangle$ and
$\langle M_{i,j} : i \leq d, j
\leq j_i \rangle ,$ satisfying the
following conditions.
\begin{enumerate}
\item If $A \subseteq M_i, \speciallns{A}{B}, B \in \kinf,$ and $B$ is
$k'$-small, then for some $j < j_i, A = A_{i,j}$ and $B$ and $B_{i,j}$ are
isomorphic over $A.$
\item $M_0 = \emptyset.$
\item For all $i \leq d+1, \cn{M_i} < 2^{k-2}.$
\item For each $i \leq d, M_{i,0} = M_i$ and $M_{i,j_i} = M_{i+1}.$
\item For each $j<j_i,$ there are $A_{i,j}, B_{i,j}$ with
\begin{enumerate}
\item $B_{i,j}$ is $k'$-small;
\item $B_{i,j} \in \kinf$;
\item $A_{i,j} \subseteq M_i$;
\item $\speciallns{A_{i,j}}{B_{i,j}}$;
\item $B_{i,j}$ is compatible with $M_{i,j}$ and $A_{i,j}$ is the subgraph
of $M_{i,j}$ induced by $B_{i,j} \cap M_{i,j}$;
\item $M_{i,j+1} = \freejoin{M_{i,j}}{B_{i,j}}{A_{i,j}}$;
\end{enumerate}
\end{enumerate}
By Sublemma \ref{closure-sublemma}, there are only
finitely many $k'$-small $B \in \kinf.$
The existence of chains satisfying the above conditions then follows
immediately from the free amalgamation property for \speciallns{}{}\
stated in Sublemma \ref{freeamalgam-sublemma}. 

We now let $M = M_{d+1}.$ It follows immediately from the construction 
that $M$ satisfies condition (M1)
above. Thus, it only remains to show that $M$ satisfies condition (M2). In
order to do so, it suffices to verify the following claim which
supplies a 
winning strategy for the Duplicator in the $d$-move \ekpgame\ from $N$ to
$M.$ 
\begin{quote}
{\em Claim}: Suppose $A = \{a_1, \ldots, a_{k'}\} \subseteq N, A' =
\cl{k'}{\infty}{A}{N}$ 
and $f$ is an embedding of $A'$ (the subgraph of $N$ induced
by $A'$) into $M_{(d+1)-i}.$ Then the pebble position with $\alpha_r$ on $a_r$
and 
$\beta_r$ on $f(a_r),$ for $1\leq r \leq k'$ is a winning position for the
Duplicator in the $i$-move \ekpgame\ from $N$ to $M.$ 
\end{quote}

We proceed to establish the claim by induction. Given $1 \leq i \leq d,$
suppose that $A, A',f,$ and the pebble position are as 
described. It suffices to show that given any move by the Spoiler, the
Duplicator can respond with a move into $M_{(d+1) - (i-1)}$ which will allow
the 
conditions of the claim to be preserved. Suppose, without loss of generality,
that the Spoiler moves $\alpha_{k'}$ onto a vertex $a \in N.$ 
Let $A'' = \cl{k'}{\infty}{\{a_1, \ldots, a_{k'-1}\}}{N}$ and
let $A''' = \cl{k'}{\infty}{\{a_1, \ldots, a_{k'-1},a\}}{N}.$ Then, by
condition 
(N2), $A''' \in \kinf$ and \speciallns{A''}{A'''}.
Then, by condition 5 on the construction of our chains defining $M,$ there is a
$B \subseteq M_{(d+1) - (i-1)}$ and an isomorphism $f'$ from $A'''$ onto $B$
with $f'$ and $f$ having identical restrictions to $A''.$ Therefore, the
conditions of the claim will be preserved, if the Duplicator plays pebble
$\beta_{k'}$ onto $f'(a).$
\end{proof}

\begin{proof}[Proof
of Corollary \ref{finbase-cor}] 
Let $k \geq 7.$ 1. Suppose,
for 
{\em reductio}, that \univk\ is finitely based with ``basis'' $\{ A_1, \ldots
, A_n \}.$ Let $k'$ be the maximum of the cardinalities of the $A_i.$ Then,
there is a sentence of \ffoe{k'} which defines \univk, contradicting Theorem
\ref{main-thm}. 

2. Suppose for {\em reductio} that the class of
$k$-extendible structures is finitely based and choose $k'$ as above with
respect to a ``basis'' for this class. As in the proof of Lemma
\ref{central-lemma}, there is a $k$-extendible graph $N$ such that each
\ffoe{k'} sentence true in $N$ has a model which is not $k$-universal and hence
not $k$-extendible. This implies that every submodel of $N$ of size at most
$k'$ is not $k$-extendible, which yields the desired contradiction. 
\end{proof}

\begin{proof}[Proof
of Theorem \ref{univinlfp-thm}]
We show that the complement of \univk\ is definable in 
least fixed point logic, which is
sufficient since the language is closed under negation.  
In fact, it is defined by a purely universal sentence.
The main idea
is to show that for all $A, A \not\in \univk$ iff either $card(A) < k-1$
or for all 
proper $k-1$-tuples $\overline{a} \subseteq A$, $\overline{a}$ is not 
$k,m$-universal for some $m \in \omega$. 
Equivalently, every proper $k-1$-tuple has finite rank.
This follows easily from the
following sequence of observations.

\begin{enumerate}
\item  For all $A$,
$A$ is $k$-universal iff there is a proper $k-1$-tuple $\overline{a}
\subseteq A$ such that $\overline{a}$ is $k$-universal in $A$.

\item  For all $A$, and every proper $k-1$-tuple $\overline{a} 
\subseteq A$, $\overline{a}$ is $k$-universal in $A$ iff
$\overline{a}$ is $k, m$-universal in $A$, for all $m \in \omega$.

\item  For every $A$ and proper $k-1$-tuple $\overline{a}$, 
if $\overline{a}$ 
has rank $m+1$ in $A,$ then there is some
set $S \subseteq \{1, \ldots , k-1 \}$ and 
formula $\varphi(x_1, \ldots , x_k) = 
\bigwedge_{i < k} x_i \neq x_k \wedge \bigwedge_{i \in S}
E(x_i, x_k) \wedge \bigwedge_{i \not\in S} \neg E(x_i, x_k)$,
such that for all $a' \in A$, if $A \models \varphi(\overline{a}a')$,
then $\overline{a}a'$ has rank $\leq m$.
\end{enumerate}

Observations 1 and 2  essentially follow immediately from the definitions.
Observation 3 may be verified by considering the $k$-extendible models.

The above conditions yield an easy inductive definition of all
the proper $k-1$-tuples that are not $k$-universal.  Call a formula
of the form of $\varphi$ above a {\em $k$-extension formula}.
Let $\varphi_1, \ldots , \varphi_t$ be the set of $k$-extension
formulas.  
By observation 3, a proper $k-1$-tuple $\overline{a}$ has 
rank 0 iff there is some $k$-extension formula $\varphi$ such that 
there is no $a'$ such that $A \models \varphi(\overline{a}a')$;
and $\overline{a}$ has rank $\leq m+1$
iff there is some $k$-extension formula $\varphi$
such that for all $a'$,
if $A \models \varphi(\overline{a}a')$, then $\overline{a}a'$
has rank $\leq m$.  

We now show how to express this definition by a least fixed point formula. Let 
$\theta(x_1, \ldots , x_{k-1})$ be the following formula:
$$
\bigvee_{i < j \leq k-1} x_i =x_j \orr 
\bigvee_{s \leq t}\forall x_k (\neg 
\varphi_s(\overline{x}x_k) \orr \bigvee_{j \leq k} R(x_1, \ldots , x_{j-1},
x_{j+1}, \ldots , x_k)).$$ 
$R$ appears positively in the 
formula, so that $\theta$ defines an inductive operator on each
graph $G, \Theta_G(X)$, that maps $k-1$-ary relations $P$ to 
$k-1$-ary relations $\Theta_G(P)$.  Let $\Theta^0_G = \Theta_G 
(\emptyset)$, and let $\Theta^{n+1}_G = \Theta_G(\Theta^n_G)$.
If $\Theta^{n+1}_G = \Theta^n_G$, then $\Theta^n_G$ is a fixed point of 
the operator.  In fact, it is the least fixed point, which we 
denote $\Theta^{\infty}_G$.  Observe that for all proper $k-1$-tuples
$\overline{a}, \overline{a} \in \Theta^{n+1}_G -\Theta^n_G$
iff the rank of $\overline{a}$ is $n$.
By the above observation, $G$ is $k$-universal iff $\Theta^{\infty}_G
=A^{k-1}$.  Therefore, the following formula defines the 
class of graphs that are not in \univk.

$$\forall x_1 \ldots x_{k-1} \bigvee_{i < j \leq k-1} x_i =x_j
\orr \forall x_1 \ldots x_{k-1} \Theta^{\infty}_G(x_1, \ldots , x_{k-1})$$

This completes the proof.
\end{proof}

\begin{proof}[Proof
of Theorem \ref{univnotextend-thm}] 
Let $k \geq 4.$ We construct
$G$ 
as follows. Let $V$ be the set of 
binary sequences of length $k,$ that is, $V$ is the set of $0,1$-valued
functions with domain $\{1, \ldots, k \}.$ For each $1 \leq i \leq k,$ let 
$V_i = V \times \{i\}$ and let $U  = \bigcup_{1 \leq i \leq k} V_i.$
$U$ 
is the set of vertices of the graph $G.$ The edge relation $E$ of $G$ is
defined 
as follows:
$$E((f,i),(g,j)) \longleftrightarrow (i \neq j \wedge f(j) = g(i)).$$
We proceed to verify that $G$ satisfies the conditions of the theorem.

First we show that $G$ is $k$--universal. Let $H$ be an arbitrary graph. We
describe a winning strategy for
the Duplicator in the \ekgame\ from $H$ to $G.$ 
At each round the Duplicator plays so as to pebble at
most one element of each $V_i.$ We may suppose without loss of generality that
all $k$ pebbles are on the board at round $s,$ that the Duplicator has played
$\beta_i$ on an element of $V_i,$ and that the map from the elements
pebbled in $H$ to the corresponding elements pebbled in $G$ is a partial
isomorphism. Suppose the Spoiler plays $\alpha_j$ onto an element $b
\in H$ at round $s+1$
and let $X$ be the set of $i$ such that there is an edge between $b$ and
the vertex of $H$ pebbled by $\alpha_i.$ Let $(f_i,i)$ be the
vertex of $G$ pebbled by $\beta_i$ at round $s.$ We must
show that the Duplicator may play $\beta_j$ at round $s+1$ onto a
vertex $(g,j) \in V_j$ such that for all $1 \leq i \leq k,$
$$E((g,j),(f_i,i)) \longleftrightarrow i \in X.$$
It is clear that $(g,j)$ satisfies this condition when $g$ is defined as
follows: $g(i) = f_i(j),$ if $i \in X$; $g(i) = 1 - f_i(j),$ if $i \not\in X.$
This completes the proof that $G$ is $k$-universal. 

Let $H \subseteq G$, and suppose, for {\em reductio}, that $H$ is
$k$-extendible.  It is easy to verify that any graph $H$ is
$k$-extendible iff
for all $j$-tuples $\overline{a}$ in $H$, $j \leq k$,
$\overline{a}$ is $k$-universal in $H$\@.  To establish the 
contradiction, we show that there are $a_1, a_2 \in H$ such 
that $(a_1, a_2)$ is not $k$-universal in $G$, which 
immediately implies that $(a_1,a_2)$ is not $k$-universal in $H$ either.

The cardinality of any $k$-extendible graph is $\geq k+1$, so there
is an $l \leq k$ such that $H$ contains two vertices,
$(f_1, l), (f_2, l)$, in $V_l$.  Let $w' = \{j \mid
j \neq l \mbox{ and } f_1(j) \neq f_2(j)\}$
and let $w'' = \{j \mid j \neq l \mbox{ and } f_1(j) = f_2(j)\}$.
Let $w =w'$, if $|w'| \leq |w''|$, and let $w = w''$, otherwise.
Observe that $|w| \leq (k-1)/2,$ which is $< k-2$ for all $k \geq 4.$
We now show that $(f_1, l), (f_2, l)$ is not $k, |w| +1$-universal
in $G$\@.  Suppose that $w=w'$.  
Let $\theta(x_1, \ldots , x_{|w| +3}) = $
$$\bigwedge_{1 \leq i<j\leq |w|+3} x_i \neq x_j
\wedge \bigwedge_{3\leq i \leq |w|+3} (E(x_1, x_i) \wedge \neg E(x_2, x_i))
\wedge \bigwedge_{3\leq i < j \leq |w| + 3}E(x_i,x_j).$$  
(Note that $|w|+3 \leq k$, since $k\geq 4$.) Observe that 
for any $|w|+3$-tuple $\overline{a} = (a_1, \ldots , a_{|w|+3})$
such that $a_1 = (f_1, l)$ and $a_2 = (f_2, l)$,
$G \not\models \theta(\overline{a})$.  If we let $$\varphi(x_1,x_2)
= \exists x_3 \ldots x_{|w|+3}\theta(x_1, \ldots , x_{|w|+3}),$$
then it follows that $G \not\models \varphi((f_1,l),(f_2,l))$.
Therefore $((f_1,l),(f_2,l))$ is not $k,|w|+1$-universal in $G$\@.
The argument for $w = w''$ is similar.
\end{proof}

The above construction may be extended to arbitrary finite relational
signatures.

\begin{proof}[Proof
of Theorem \ref{kineqminkuniv-thm}]
Let $k \geq 6.$  For all $n \geq 4k$, we construct
graphs $G_n$ such that:
\begin{enumerate}
\item  $G_n$ is $k$-universal.
\item  For all $H \subseteq G_n$, if $H$ is $k$-universal, then 
  the diameter of $H$ is $\geq\lfloor (n-1)/2\rfloor /(k-1)$.
\end{enumerate}
(Recall that the diameter of a graph is the maximum distance 
between any two vertices if it is connected, and $\omega$ otherwise.
It is an easy exercise to show that for $k \geq 3$, every minimal
$k$-universal graph is connected.)
This immediately yields the fact that there are 
minimal $k$-universal models of arbitarily large finite diameter.
It is easy to check that the property of having finite diameter 
$=d$ is expressible in $L^3$, which implies that any two graphs 
with different diameters are $L^k$-inequivalent.

The graphs $G_n$ are based on a modification of the construction
from the proof of Theorem 3.  Let $V$ be the set of functions from 
the interval $\{-(k-2), \ldots , k-2\}$ into $\{0, 1\}$.
For each $m, 0 \leq m \leq n-1$, let $V_{m} = \{0, 1\} 
\times V \times \{ m\}$.
The set of vertices of $G_n$ is $\bigcup_m V_{m}$.  The edge 
relation on $G_n$ is defined as follows.  For all 
$m,m'$, $a \in V_m, a' \in V_{m'}$, if $m= m'$ or 
$k \leq m-m' \leq n-k(mod\ n)$,
then $\neg E(a,a').$  If $0 < m-m' <k-1 (mod\ n)$,
and $a = (\delta , f, m), a' =(\delta ', f', m')$, with $\delta , 
\delta ' \in \{0, 1\}$ and $f,f' \in V$, then $E(a,a')$
iff $f'(m-m')=f(m'-m)$. (Here, subtraction is $modulo\ n$.)
Finally, if $m-m' = k-1 (mod\ n)$, then 
$E(a,a')$ iff $\delta = 1$.  In this case, each $a \in V_m$
is either adjacent to every vertex in $V_{m'}$ or to none of them.
If $m' = m + \lfloor (n-1)/2\rfloor$, then the distance $d(a,a') \geq
\lfloor(n-1)/2\rfloor /(k-1)$.
Observe also that for all $l \leq n-1$,
there is an automorphism of $G_n$ taking each $V_m$
to $V_{m+l}$.  (All indices are $modulo\ n$.)

First we show that $G_n$ is $k$-universal.  Let $G'$ be an arbitary
graph.  It suffices to prove that the D wins the $\exists^k$-game
from $G'$ to $G_n$.
By an argument similar to the one given in the proof of
Theorem 3, it is easy to see that the D can play so that in each
round $i \leq k$, she plays a pebble on a vertex in $V_i $.
We now argue by induction that in each subsequent round $j > k$,
she can maintain the following condition:  there is some $l \leq n$
such that there is exactly one pebble on each $V_m$, for $m$ such that
$0 \leq m-l \leq k-1(mod\ n)$.  The basis step is already taken care of.
Suppose that in round $j$, the D has a single pebble in each 
vertex set $V_l,  \ldots , V_{l+(k-1)}$.  We consider 
two cases.  One, the S replays the pebble $\alpha_i$ whose pair 
$\beta_i$ in $G_n$
is on an element of $V_l$.  It is easy to verify that the
D can respond by playing $\beta_i$ on a vertex in $V_{l+k}$.
Observe that the D's pebbles are now on $V_{l+1},  \ldots ,
V_{l+k}$, as desired.
Two, the S replays any other pebble $\alpha_{i'}$, whose pair 
$\beta_{i'}$ is on some element
of $V_{l'}, l \neq l'$.  The D can  respond by replaying 
the pebble on some other element of $V_{l'}$.  Again, 
that this is possible essentially follows from the proof of Theorem 3. 

Next we argue that any $k$-universal $H \subseteq G_n$ has diameter
$\geq \lfloor (n-1)/2\rfloor /(k-1)$.  In particular, it is 
sufficient to prove $H$ must contain a vertex from each $V_m, m\leq
n-1$.  Let $A$ be any $k$-extendible graph.  The argument proceeds
by establishing that, in the $\exists^k$-game from $A$ to $H$,
the S can eventually force the D to play a pebble on a vertex in 
each $V_m \cap H$.  If $V_m \cap H = \emptyset$, for some $m$,
then the D loses.

In rounds 1 through $k$, the S plays on a $k$-clique in $A$.  
For every $k$-clique in $G_n$, and hence also in $H$, there is an
$m \leq n-1$ such that each $V_{m'}, 0 \leq m'-m \leq k-1(mod\ n)$, contains
exactly one element from the clique.  Therefore, after $k$ rounds,
the D must have a single pebble on each of $V_m, \dots ,
V_{m+(k-1)}$, for some $m$.  It suffices to show that the S
can force the D to play so that exactly one pebble occupies a 
vertex in each set $V_{m+1},
\ldots , V_{m+k}$, since by iterating this strategy, he can
force the D to play onto each $V_{l}$. 

To simplify the notation, we assume $m=0$ and that each pebble $\beta_i,
0 \leq i \leq k-1$, is on a vertex in $V_{i}$.  Let 
$b_i = (\delta_i, f_i, i), \delta_i \in \{0, 1\}, f_i \in V$,
be the element pebbled by $\beta_i$. 
In round $k+1$, the S replays 
pebble $\alpha_0$ and places it on an element $a \in A$ such 
that $E(a, \alpha_1)$ and 
for $i \in \{2, \ldots , k-1\}$, $E(a, \alpha_i)$ iff
$\delta_i = 0$. (Here we abuse notation
and use
$\alpha_j$ to refer also to the element on which the pebble is 
located.)  
Since $\alpha_0$ and $\alpha_1$ are now adjacent in $A$, the D has to 
play $\beta_0$ on some element in a set $V_l$, for
$-(k-2) \leq l \leq k(mod\ n)$, so that it is adjacent to $\beta_1$.

By the condition that for $i \in \{2, \ldots ,
k-1\}$, $E(a, \alpha_i)$ iff $\delta_i = 0$, the D cannot play
in $V_l$, for $-(k-3) \leq l \leq 0(mod\ n)$.  
If the D plays the pebble in $V_k$, then 
the S has succeeded. Suppose that the D plays
$\beta_0$ on an element of $V_{-(k-2)}$.  We now claim that 
there is no 3-clique in $G_n$ [$H$] each of whose elements is 
adjacent to both  $\beta_{k-1}$ and $\beta_0$.
This is because $(i)$ the only elements of $G_n$ that are adjacent
to vertices in both $V_{-(k-2)}$ and $V_{k-1}$ are members 
of either $V_{0}$ or $V_{1}$, and  $(ii)$ there is no 3-clique
in $V_{0} \cup V_{1}$.  Thus the S can force a win in 3 moves
by replaying pebbles $\alpha_1, \alpha_2, \alpha_3$ so that they
occupy a 3-clique each of whose elements are adjacent to $\alpha_0$
and $\alpha_{k-1}$.  

The remaining case occurs when the D plays $\beta_0$ on a vertex
in $V_{j}$, for $1 \leq j \leq k-1$.  Without loss of generality,
let $j = k-2$, and let $b'$ be the vertex now occupied by 
$\beta_0$.  Let $w'= \{i \mid 1 \leq i \leq 3 \mbox{ and }
E(b_{k-2}, b_i) \mbox{ iff } E(b', b_i)\}$ and 
$w''= \{i \mid 1 \leq i \leq 3 \mbox{ and }
E(b_{k-2}, b_i) \mbox{ iff } \neg E(b', b_i)\}$.  
Again without loss of generality, suppose that 
$|w'| \geq 2$ and $w' =\{1, 2\}$. 
By exploiting the fact that $\beta_0$ and $\beta_{k-2}$ both
occupy vertices in $V_{k-2}$, the S can now force
the D to play $\beta_2$ onto $V_k$.

The S first places $\alpha_2$
on a vertex such that for all $j, 1 \leq j \leq k-1, j\neq 2$,
$E(\alpha_2, \alpha_j)$, and $\neg E(\alpha_2, \alpha_0)$.
It is easy to see that the D must put $\beta_2$ on either $V_0$ or
$V_k$.  
In the first case, the S responds by playing
$\alpha_1$ so that for all $j, 2 \leq j
\leq k-1$, $E(\alpha_1, \alpha_j)$ and $\neg E(\alpha_1, \alpha_0)$.
The D now loses immediately.  The only vertices adjacent to each
$\beta_j, 2 \leq j \leq k-1$, are elements of $V_1$ or $V_2$,
but for each $b \in V_1$ or $V_2,$ $E(b, \beta_{k-2})$ iff
$E(b, \beta_0)$.
In the second case, the S then plays 
$\alpha_0$ onto a vertex such that for all $j, 1 \leq j \leq k-1$,
$E(\alpha_0, \alpha_j)$.  This compels the D to play $\beta_0$ in
$V_{2}$, so that there is a now a single pebble in each
$V_1, \ldots , V_k$,
as desired.    
\end{proof}
\bibliographystyle{amsalpha}
\providecommand{\bysame}{\leavevmode\hbox to3em{\hrulefill}\thinspace}

\end{document}